\input amstex.tex
\documentstyle{amsppt}
%\magnification\magstep1

\topmatter \title Gerbes on quantum groups \endtitle
\author Jouko Mickelsson \endauthor 
\affil Mathematical Physics, KTH, SE-106 91 Stockholm, Sweden
\endaffil 
\date August 25, 2003 \enddate 
\endtopmatter

\advance\vsize -2cm
\magnification=1200
\hfuzz=30pt
\document

\baselineskip=16pt
\redefine\l{\lambda}

\bf Abstract \rm  We discuss an approach to quantum gerbes over
quantum groups in terms of q-deformation of transition functions 
for a loop group bundle. The case of the quantum group 
$SU_q(2)$ is treated in some detail.   

\vskip 0.3in 
\bf 0. Introduction \rm 

\vskip 0.2in
In view of the recent extensive activity in the theory of gerbes and
their applications in quantum field theory it  is not surprising that 
the  question arises whether there is some sort of object which could
be  called a 'quantum gerbe'. One should understand 'quantum' here
really  meaning a deformation depending on a real parameter $q$ since 
a gerbe  already is in a sense a quantum object: It has a well defined 
'quantum  number' given by its Dixmier-Douady class  and this class in 
the  physical applications is related to the chiral anomaly in quantum 
field  theory. 

A proposal  for a quantum gerbe was discussed in [ABJS] in terms of
deformation  quantization. The underlying base space of the gerbe was 
replaced by  noncommutative space defined as a star product algebra. 
A noncommutative  version of a line bundle was defined, [JSW], and then a gerbe 
was defined as a  system of local noncommutative line bundles obeying 
a certain cocycle  condition, imitating the corresponding cocycle 
condition for the  undeformed case. 

In this paper  we follow a different route to quantum gerbes. We start 
from the  alternative definition of a gerbe (without connection and
curvature)  as a principal $PU(H)$ bundle over a manifold $M,$ where 
$PU(H)$ is  the projective unitary group of a complex Hilbert space
$H.$  The gerbe class can be  nontorsion only if $H$ is
infinite-dimensional.   
The  equivalence classes of such a bundles are parametrized by elements
in $H^3(M,\Bbb Z).$  The characteristic class is called the
Dixmier-Douady  class. 

Concretely, the principal $PU(H)$ bundle can be obtained from a loop 
group bundle using a projective unitary representation of the loop
group. In this paper we shall construct 'quantum loop group bundles' 
using transition functions in a group of matrices with entries in a 
quantum group.  

The  local approach to gerbes  consists of an open cover $\{U_{\l}\}$ of $M$ and 
a set of line bundles $L_{\l\l'}$ over the intersections $U_{\l}\cap 
U_{\l'}$ with given isomorphisms 
$$L_{\l\l'}\otimes L_{\l'\l''}= L_{\l\l''}\tag1$$ 
on triple intersections. The relation to the global description is 
as follows. There is a canonical central extension 
$$1\to S^1 \to U(H) \to PU(H) \to 1$$ 
of the projective unitary group which gives a canonical complex line 
bundle $L$ over $PU(H).$ The transition functions $\phi_{\l\l'}: 
U_{\l}\cap U_{\l'} \to PU(H)$ of the principal bundle can then be used 
to pull back the line bundle $L$ over $PU(H)$ to a local line bundle 
$L_{\l\l'}$ over $U_{\l}\cap U_{\l'}.$ The group structure on $U(H)$ 
defines an identification of $L_{\l\l'}\otimes L_{\l'\l''}$ with 
$L_{\l\l''}.$ 

In the case of the quantum group $SU_q(2)$ we shall see how the 
the local line bundle approach is related to the loop group bundle 
construction. 

\vskip 0.3in
\bf 1.  The gerbe over $SU(n)$: The undeformed case  \rm 

\vskip 0.2in
Let  us next consider the case when $M=G$ is a compact connected Lie group. 
Let $P$  be the space of all smooth paths $f$ in $G$ starting from
$f(0)=1$  and with an arbitrary endpoint $f(1)\in G.$ We also require
that $f^{-1}df$ is a smooth periodic function, that is, it defines 
a vector potential on the unit circle $S^1.$ Thus we may identify $P$
as the contractible space $\Cal A$ of smooth vector potentials on
$S^1$ with values in the Lie algebra of $G.$

$P$ is the total space of a principal bundle over $G$ with fiber 
equal to the group $\Omega G$ of smooth based loops in $G. $ Based 
means that $f(0)=f(1)=1.$ Since $P=\Cal A$ is contractible, it is 
a universal bundle for $\Omega G.$ 
Assume that $\psi: \Omega G \to PU(H) $ 
is a projective representation of $\Omega G$ in the Hilbert space $H.$ 
We can then define an associated principal $PU(H)$ bundle over $G$ 
in the usual way, with total space $Q= P\times_{\psi} PU(H).$ 
This is the way gerbes appear in canonical quantization in field
theory, [M, CM1-2, CMM].

We shall \it define  \rm a quantum gerbe over a quantum group $SU_q(N)$ 
in terms of 'local transition function' with values in a loop group. 
This means that the transition functions are loops in a group of 
matrices with entries in the Hopf algebra $SU_q(n).$ 

To warm up, we start from the classical case by giving explicit 
formulas for the transition functions. In the case of $SU(n)$ it is 
sufficient to select $n$ open sets to cover the base. We choose 
$n$ different points $\l_1,\dots \l_n\neq 1$ in the unit circle in the 
complex plane, ordered counter clockwise, such that the product 
$\l_1\l_2\dots\l_n\neq 1.$   For $i=1,2,\dots n$ let 
$U_i$ be the subset of $SU(n)$ consisting of matrices $g$ such that 
$\l_i$ is not an eigenvalue of $g.$ This gives an open cover: If 
$g$ is not in any of the sets $U_i$ then all $\l_i$'s are eigenvalues 
of $g$ and so det$(g)= \l_1\l_2\dots\l_n\neq 1,$ a contradiction.

An each open set $U_i$ we have a trivialization of the bundle $P\to 
SU(n).$ Let $h_i(t)$ be any fixed smooth contraction of the set 
$S^1\setminus \{\l_i\}$ with $h_i(0)$ the constant map sending 
$S^1\setminus\{\l_i\}$ to the point $1$ and $h_i(1)$ the identity 
map. Let $d=(d_1,\dots, d_n)$ be a diagonalization of $g\in U_i,$
$g=AdA^{-1}.$ Then we set  
$$\psi_i(g)(t)= A (h_i(t)(d_1),\dots, h_i(t)(d_n)) A^{-1}$$
and this defines a path in $U_i$ joining $g$ to the neutral element, 
i.e., we have a local section $\psi_i: U_i \to P.$ 

The transition function on $U_i\cap U_j$ is fixed by $\psi_j(g) =
\psi_i(g) \phi_{ij}(g)$ and it takes values in $\Omega G.$ 

\bf Remark 1 \rm   Allowing the paths lie in the bigger group $GL(n,\Bbb C)$ 
leads to a  technical simplification in the construction of transition 
functions.  First, for each index $i$ we have the contraction of $U_i$
to the point $-\l_i\cdot 1$ given by 
$$\psi_i(g)(t)= -(1-t)\l_i \cdot 1 + t g.$$
Note that  the matrix $\psi_i(g)(t)$ is really invertible for each 
$g\in U_i$  and $0\leq t\leq 1$ by the spectral property Spec$(g)
\subset S^1.$  Since $GL(n,\Bbb C)$ is connected, we can deform 
the constant  map $\psi_i(g)(0)= -\l_i \cdot 1$ to the constant map 
$g\mapsto 1$  in an obvious way. Putting these together we obtain 
a homotopy  connecting the constant map $g\mapsto 1$ to the indentity 
map $g\mapsto g$  in $U_i.$ The transition functions are defined as
before, $\phi_{ij}(g)= \psi_i(g)^{-1} \psi_j(g),$ now with values 
in the loop group of $GL(n,\Bbb C).$ 

A third way to define the trivializations and transition functions is
to use functional operator calculus.  If $\l\in S^1$ is not in the
spectrum of a unitary matrix $g$ then the logarithm of $g$ can be
defined as
$$\log(g)= \int_{\gamma} \frac{\log(z)}{g-z}dz,\tag2$$ 
where $\gamma$ is a closed loop in the complex plane $\Bbb
C\setminus\{0\}$ which encircles
every point in $S^1\setminus\{\lambda\}$ with winding number one; the 
loop crosses the point $\l$ twice such that the ray $R$ from the origin 
to $\lambda$ is tangential to $\gamma$ at $\l$ but otherwise $\gamma$
does not cross the ray; then we can define the branch of $\log(z)$ as
an analytic function with a cut along $R$ and $(g-z)^{-1}$ is
nonsingular along $\gamma.$ Denoting $X=\log(g)$ we have a smooth path 
$ \psi_i(g)(t)= e^{tX}$ connecting $1$ to the point $g\in U_i$ for
$\l= \l_i.$

As mentioned in the beginning, the local line bundle construction of 
a gerbe over $SU(n)$ is obtained by pulling back the (level $k$)
central extension of the loop group by the transition functions. 
However, there is simpler construction which leads to the same gerbe. 

As shown in [CM1] the basic gerbe (level $k=1$) is constructed using 
a family of Dirac operators parametrized by points on $G.$ For 
$G=SU(n)$ the Dirac operator $D_g$ attached to $g\in G$ has domain 
consisting of smooth functions on the interval $[0,1]$  with values 
in $\Bbb C^n$ and boundary conditions $\psi(1)= g\psi(0).$ The Dirac 
operator is simply the differentiation $D_g= -i\frac{d}{dx}$ on the 
interval.  If $g=1$ the spectrum has multiplicity $n$ and consists of
the numbers $2\pi m, m\in \Bbb Z.$ For arbitrary $g$ the spectrum 
consists of $2\pi m + \mu_i$ where $\mu_1, \dots , \mu_n$ are the 
eigenvalues of $-i\log(g).$ 

If now $\l,\l'$ are two distinct points on the unit circle, $\l= 
e^{i\mu}, \l'=e^{i\mu'},$ with $0\leq \mu,\mu' <2\pi,$ 
then the eigenvalues of the Dirac operator
$D_g$ with $g\in U_{\l}\cap U_{\l'}$ in the spectral interval 
$]\mu,\mu'[$ match the eigenvalues of $g$ in the segment $]\l,\l'[$ 
of the unit circle. It follows that the determinant line $DET_{\mu\mu'}$
for the Dirac operators $D_g$ is naturally isomorphic  to the top 
exterior power $\Lambda^{top}(E_{\l\l'})$ of the corresponding spectral subspace 
for $g.$  Combining with [CM1] we get the equivalences:

\proclaim{Theorem} The basic gerbe over $G=SU(n),$ corresponding to
the Dixmier-Douady class given by the generator of $H^3(SU(n),\Bbb Z),$ 
can be given in three equivalent ways: 
\roster \item as the $PU(H)$ bundle over $G$ obtained as an associated bundle 
to the universal bundle $\Cal A\to G$ through the basic representation 
of the affine Kac-Moody group based on $G$ 
\item as the local system of complex line bundles defined as above by Dirac
operators on the unit interval, parametrized by boundary conditions
$g\in G$ 
\item as the system of local line bundles formed from top exterior 
powers of the spectral subspaces of $g\in G$ corresponding to 
open segments of the unit circle. 
\endroster
\endproclaim

Let us consider the case $G=SU(2)$ in more detail. It is sufficient to 
choose a cover consisting of two open sets $U_{\pm}$ consisting of
points $g\neq \pm 1.$ The overlap is homotopic to the equator $S^2$ 
in $SU(2) = S^3.$ The basic transition function for a loop group
bundle on $SU(2)$ is then 
$$\phi(x)(t) = \cos(2\pi t) +i x \sin(2\pi t) \text{ for $0\leq
t\leq\frac12$} \tag3$$ 
and $\phi(x)(t)= diag(e^{2\pi it}, e^{-2\pi it})$ for $\frac12 \leq
t\leq 1.$ Here 
$x=\left(\matrix x_3& x_1+ix_2 \\ x_1-ix_2& -x_3\endmatrix\right).$
Actually, to make the function smooth at $t=1/2$ and the end
points one should
replace the variable $t$ by $f(t)$ where $f$ is any smooth monotonous
function on the interval $[0,1]$ such that $f(0)=0, f(1)=1,$ and all 
the derivatives vanish at the points $t=0,1/2,1.$  

\vskip 0.3in 
\bf 2. The gerbe over the quantum group $SU_q(n)$ \rm

\vskip 0.2in
We assume $0< q < \infty$  and denote by $SU_q(n)$ the standard Hopf algebra 
of the quantum special  unitary group. It is given in terms of
generators and relations  as follows. The generators are elements
$g_{ij}$ indexed by $i,j=1,2,\dots n$ with defining relations

$$ \align &g_{im}g_{ik}= q g_{ik} g_{im}, \,\, g_{jm}g_{im}=qg_{im}g_{jm}\\  
g_{im} g_{jk}& = g_{jk}g_{im}, \,\,\,
g_{ik}g_{jm}-g_{jm}g_{ik}=(q^{-1}-q)g_{im} g_{jk}\tag4\endalign$$ 
for $i<j$ and $k<m.$ In addition, the quantum determinant
$$det_q= \sum_{\sigma\in S_n} (-q)^{-\ell(\sigma)} g_{1\sigma(1)} 
g_{2\sigma(2)} \dots g_{n\sigma(n)} \tag5$$ 
is identified as the unit element in the algebra. 

 The counit is
defined by $\epsilon(g_{ij})= \delta_{ij}$ and the antipode is
defined as 
$$S(g_{ij}) = (-1)^{i+j} X_{ji},  \tag6$$ 
where $X_{ij}$ is the $(n-1)\times(n-1)$ quantum minor of the matrix 
$(g_{ij}),$ i.e., the quantum determinant of the submatrix obtained by
deleting the $i$:th row and $j:$the column from $(g_{ij}).$ The star
algebra structure is given as
$$g_{ij}^* = S(g_{ji}).\tag7$$ 
In particular, by the definition of an antipode,
$\sum_{(x)}  x'S(x'')=\epsilon(x)\cdot 1$ where $\Delta(x)= \sum_{(x)} x'\otimes
x'',$ 
and so
$$g_{ij} g_{ik}^* = \delta_{jk} \cdot 1.\tag8$$
The coproduct $\Delta$ is defined by matrix multiplication,
$$\Delta(g_{ij})= \sum_k g_{ik}\otimes g_{kj}.$$ 
For more details, see [KS].  

In the case $n=2$ the standard notation is $a=g_{11}, b=g_{12}, 
c=g_{21}, d=g_{22}$ with $b^*= -qc,  a^*=d.$ The determinant condition 
is $ad-qbc=1$ and the commutation relations are 
$$ab=qba,\,\, ac=qca, \,\,bd=qdb, \,\,cd=qdc, \,\,bc=cb, \,\,ad-da= (q-q^{-1})bc. \tag9$$ 

We shall now define a 'quantum gerbe' over the quantum group $SU_q(2)$ 
as a quantum projective bundle over $SU_q(2).$ This is achieved by
giving the transition function on the 'equator' of $SU_q(2).$ 
The quantum equator is defined as a quotient algebra of $SU_q(2),$ 
[HMS]. Let $I\subset SU_q(2)$ be the 2-sided ideal generated by 
the single element $b-b^*= b+qc.$ Then one can show that the standard 
quantum sphere $S^2_q$ is isomorphic to $SU_q(2)/I, $
[HMS]. Explicitly,
the generators of $S^2_q$ are $K=K^*$ and $L$ with the defining
relations 
$$LL^* +q^2K^2 =1, \,\,L^*L +K^2=1, \,\, LK=qKL.\tag10$$ 
One can check that $K=c+I, L=a+I$ in $SU_q(2)/I$ indeed satisfy the 
relations above. In the classical case, $q=1,$ going to the quotient 
means that we set one of the coordinates, namely Im$(b),$ equal to
zero reducing the 3-sphere to the equatorial 2-sphere. 

The transition function is defined as in the undeformed case: 
Define the quantum $2\times 2$ matrix 
$$ x= \left(\matrix qK& L\\ L^* & -K \endmatrix\right).\tag11$$
It is easy to check from the defining relations that $x^*=x$ 
(combination of matrix transposition and star operation on $S^2_q$) 
and that $x^2=1.$ It follows that $\cos(2\pi t) +i x \sin (2\pi t)$ 
($0\leq t\leq 1/2$) is unitary and in combination with the path $t\mapsto diag(e^{2\pi
it}, e^{-2\pi it})$ ($1/2\leq t\leq 1$) 
defines a loop in unitary $2\times 2$ matrices  with coefficients in the algebra $S^2_q.$ 

In the classical undeformed case the transition function can be
extended from the equator to the open set 
$U_{+-}=SU(2)\setminus \{\pm \sigma \}$ where
$\sigma =\left(\matrix\ 0& i\\ i& 0\endmatrix\right).$ 
First one extends the function $x$ to $U_{+-}$ by writing 
$$x= f^{-1}\left(\matrix \frac12(b+b^*) & a^*
\\ a & -\frac12(b+b^*)\endmatrix\right),\tag12$$
where $f=[1+\frac12(b-b^*)^2]^{1/2}.$ This function is singular at 
the points where Im$(b) = \pm i,$ that is, at opposite poles of 
$S^3=SU(2).$  It satisfies $x^*=x= x^2$ and its restriction to the 
equator defined by Im$(b)=0$ is equal to $\left(\matrix b& d \\ a&
-b \endmatrix\right).$ We may also view $x$ as a matrix with entries in an algebraic 
extension of the commutative algebra of functions on $SU(2)$ by a 
single element $f$ satisfying the defining relation $f^2=
1+\frac14(b-b^*)^2.$ That is, we add certain functions with
singularities at the poles to the algebra generated by the elements 
$a,b,c,d.$

The generalization of  the above formula to the q-deformed case is
given as
$$ x= \left(\matrix \frac{1}{2q} (b+b^*)f^{-1}_q & a^* f^{-1} \\ 
f^{-1} a & -\frac12 (b+b^*) f^{-1}  \endmatrix\right) \tag13$$ 
where $f$ is same as before and $f_q= [1+
\frac{1}{4q^2}(b-b^*)^2]^{1/2}.$ This means that the matrix $x$ has
entries in an algebraic extension of $SU_q(2)$ by the elements $f,f_q$
satisfying the relations 
$f^2= 1+\frac14 (b-b^*)^2$ and likewise for $f_q.$ 
 
As in the undeformed case, we can think of the quantum gerbe coming 
from a K-theory class defined by a family of Dirac type operators. 
The family of operators is again given by the differentiation 
$D= -i\frac{d}{dt}$ but now acting on 2-component spinors with 
coefficients in the quantum group  $SU_q(2).$ The family of boundary 
conditions is
$$\psi(1)= g\psi(0)$$ 
where $g$ is the quantum matrix $g=\left(\matrix a&b\\c&d\endmatrix\right).$   
So in a sense the K-theory class in $K^1(SU_q(2))$ is tautological: It
is the unitary $2\times 2$ matrix $g$ given by the 'coordinates' of
the quantum group. Note that the K-theory of $SU_q(2)$ when $q>0$ is
the same as for $q=1,$ $K^1(SU_q(2))=\Bbb Z = K^0(SU_q(2)),$ [MNW2]. 

The quantum line bundle on the equator $S^2_q$ is defined as the 
projection $P=\frac12(1+x)$ with coefficients in $S^2_q.$ The 
property $P^2=P$ follows from $x^2=1.$ Of course the case of $SU_q(2)$ 
is very simple since there is only one 'overlap' $S^2_q$ and we do not
need to bother about the generalization of (1) to the quantum group
case. 

\bf Remark 2 \rm The rank one projectors in the $q$ deformed case define 
only right $S^2_q$ modules and not bimodules since the left
multiplication by the elements in the algebra $S^2_q$ does not commute
with the multiplication by $P$ on $S^2_q\oplus S^2_q.$ For this reason 
the tensor products of line bundles over $S^2_q$ are not canonically
defined.  
 
Since there is no underlying smooth manifold, the Dixmier-Douady 
class cannot be defined in terms of de Rham forms as in the undeformed 
case. Instead, one can use a cyclic cocycle $c_3$ to compute the 
quantum invariant of the gerbe by pairing $c_3$ with the K-theory 
class of the unitary matrix $g.$ This was in fact already done in 
[Co] (in the case $0<q<1$) and I will not repeat (the rather complicated) calculations
here. 

Next we want to generalize the construction to $G=SU_q(n).$ 
This is achieved by a generalization of the method in Remark 1. 
Let $\l\in S^1.$ Define an algebraic extension $SU_q(n)_{\l}$ of
$SU_q(n)$ by requiring that the equation $(g-\l)h=1$ has a solution 
$h$ as a $n\times n$ matrix with entries in the extended algebra. 
In the classical case this means that we allow singularities for 
det$(g-\l),$ i.e., we restrict the domain of functions on $SU(n)$ to 
the open subset specified by $\l \notin Spec(g).$ 

Actually, we need to consider the algebra $SU_q(n)$ as a $*$-algebra 
completion of the algebra defined by the relations (4). Then $g-t\l$ 
has an inverse for any real number $t\neq 1$ since $g$ is a unitary 
matrix. We have then a homotopy 
$$ \psi_t(g,\l)= -(1-t)\l + tg\tag14$$ 
connecting $\psi_0(g,\l) = -\l$ to the identity map $g\mapsto
g=\psi_1(g,\l).$ The homotopy is defined in the space of nonsingular
matrices with entries in $SU_q(n)_{\l}.$ 
We extend this homotopy to a path connecting the constant map 
$g\mapsto 1$ to the identity $g\mapsto g$ in an obvious way, 
along the path $\psi_t= (-\l)^t,$ $0\leq t\leq 1$ connecting the neutral 
element to $\l.$  

The 'transition functions' are defined as in the undeformed case, 
$\phi_{\l\l'}(g)(t)= \psi_t(g,\l)\psi_t(g,\l')^{-1}.$
The transition functions are loops of $n\times n$ matrices with
coefficients in the extension $SU_q(n)_{\l\l'}$ defined by the
inverses $(g-\l)^{-1}, (g-\l')^{-1}.$

\vskip 0.2in
\bf Acknowledgements \rm  This work was supported by Erwin Schr\"odinger Institute for 
Mathematical Physics and CPT in Luminy (Universite de Provence) and  I want 
to thank especially Thomas Sch\"ucker for hospitality in Luminy.

\vskip 0.3in
\bf References \rm 

[ABJS]  P. Aschieri, I. Bakovic, B. Jurco, P. Schup:
Noncommutative gerbes and deformation quantization.
hep-th/0206101 

[CM2] A. L. Carey and J. Mickelsson:  The universal gerbe, Dixmier-Douady class, and gauge theory.
     Lett.Math.Phys. \bf 59, \rm 47-60 (2002).  hep-th/0107207 

 [CM1] A. L. Carey and J. Mickelsson: A gerbe obstruction to
quantization of fermions on odd dimensional manifolds with boundary. 
Lett.Math.Phys.\bf 51, \rm 145-160,(2000).  
hep-th/9912003 

[CMM] A.L.  Carey, J. Mickelsson, and M. Murray:  Bundle gerbes applied to
field theory. Rev.Math.Phys.\bf 12, \rm 65-90 (2000).  hep-th/9711133 

[Co] A. Connes:  Cyclic Cohomology, Quantum group Symmetries and the Local Index Formula for $SU_q(2).$
    math.QA/0209142

[HMS] P.M. Hajac, R. Matthes, W. Szymanski:  Graph C*-algebras and Z/2Z-quotients of quantum spheres.
          math.QA/0209268   

[JSW] B. Jurco, P. Schupp, J. Wess: Noncommutative line
bundles and Morita equivalence.  
Lett.Math.Phys.\bf 61, \rm 171-186 (2002). 
hep-th/0106110 

[KS] L.I. Korogodski and Ya. S. Soibelman: \it  Algebras of functions
  on quantum groups. Part I. \rm Mathematical Surveys and Monographs,
  \bf 56. \rm  American Mathematical Society, Providence, RI, 1998.

[MNW1] T. Masuda, Y. Nakagami, and J. Watanabe: Noncommutative
differential geometry on the quantum two sphere of Podles. I: An
algebraic viewpoint. K-Theory \bf 5, \rm 151-175 (1991) 

[MNW2] T. Masuda, Y. Nakagami, and J. Watanabe: Noncommutative
differential geometry on the quantum $SU(2).$ I: An algebraic
viewpoint. K-Theory \bf 4, \rm 157-180 (1990) 

[Mi]  J. Mickelsson:  On the Hamiltonian approach to commutator
anomalies in $(3+1)$ dimensions. Phys.Lett. \bf B241, \rm 70-76, (1990)

 \enddocument   
 
\bf 10. The quantum principal bundle: Case of $SU(2)$ bundle over
$S^4.$ \rm 

\vskip 0.2in
In the classical undeformed case the space of functions on a principal 
bundle $P$ over $S^4$ with the structure group $SU(2)$ can be
described as follows. Let $D_{\pm}$ be a pair of 4-disks with a common 
boundary $S^3.$ The principal bundle $P$ is obtained by glueing
together $P_{\pm} = SU(2) \times D_{\pm}$ along the boundary, using 
the identification $(g, g') \sim (gg',g')$ on $SU(2) \times S^3,$
where we have used the identification $SU(2)=S^3.$ (Here we consider
the case of winding number $=1.$) 

We consider the sphere $S^4$ as a suspension of $S^3.$ Concretely, 
we define the commutative algebra consisting of polynomials
in  a real variable $t$ 
and a set of new variables $g_{ij}$ (which should be thought of as a 
scaling factor times the old group variables $g_{ij}),$ such that 
$$g_{ik} g_{jk}^* = (1-t^2)\delta_{ij}.$$ 
The points $t=\pm 1$ correspond to the opposite poles of the sphere 
$S^4$ whereas $t=0$ is the equator $S^3.$  

The total space $P$ is defined algebraically as the algebra of
polynomial functions in variables $t$ and $a_{ij}, b_{ij}$ with 
the defining relations
$$a_{ij}a_{ik}^*=(1-t^2)\delta_{jk}, \,\,   b_{ij}b_{ik}^*= t^2\delta_{jk}.$$ 

The right action of $SU(2)$ is then given by the coproduct on
functions, 
$$\Delta(a_{ij})= a_{ik}\otimes g_{kj}, \,\, \Delta(b_{ij})=
b_{ik}\otimes g_{kj}$$ 
and $\Delta(t)=t\otimes 1.$ The subalgebra $B$ of $P$ generated by the elements
$$h_{ij}= a_{ik} b_{jk}^*$$
and by $t$ is invariant with respect to the coaction of $SU(2)$ and is 
isomorphic to the algebra of polynomial functions on $S^4.$ 

\vskip 0.3in

\bf 11. The quantum $SU_q(n)$ bundle over $\Sigma SU_q(n)$ \rm 

\vskip 0.2in
The quantum group $SU_q(n)$ is defined as the polynomial algebra in
the  generators $g_{ij}$ with $i,i=1,2,\dots,n$ and the defining
relations
$$ \align &g_{im}g_{ik}= q g_{ik} g_{im}, \,\, g_{jm}g_{im}=qg_{im}g_{jm}\\  
g_{im} g_{jk}& = g_{jk}g_{im}\,\,
g_{ik}g_{jm}-g_{jm}g_{ik}=(q^{-1}-q)g_{im} g_{jk}\tag2.1\endalign$$ 
for $i<j$ and $k<m.$ In addition, the quantum determinant
$$det_q= \sum_{\sigma\in S_n} (-q)^{-\ell(\sigma)} g_{1\sigma(1)} 
g_{2\sigma(2)} \dots g_{n\sigma(n)} \tag2.2$$ 
is identified as the unit element in the algebra. The counit is
defined by $\epsilon(g_{ij}= \delta_{ij}$ and the antipode is
defined as 
$$S(g_{ij}) = (-1)^{i+j} X_{ji},  \tag2.3$$ 
where $X_{ij}$ is the $(n-1)\times(n-1)$ quantum minor of the matrix 
$(g_{ij}),$ i.e., the quantum determinant of the submatrix obtained by
deletetin the $i$:th row and $j:$the column from $(g_{ij}).$ The star structure is given as
$$g_{ij}* = S(g_{ji}).\tag2.4$$ 
In particular, by the definition of an antipode,
$xS(x)=\epsilon(x)\cdot 1$
and so
$$g_{ij} g_{ik}^* = \delta_{jk} \cdot 1.\tag2.5$$
The coproduct is defined by matrix multiplication,
$$\Delta(g_{ij})= \sum_k g_{ik}\otimes g_{kj}.$$ 
For more details, see Soibelman:..... 

The suspension algebra $B=\Sigma SU_q(n)$ is generated by $g_{ij}$ and
$s$ with the modified relations 
$$g_{ij} g_{ik}^*= (1-s^2) \delta_{jk}.\tag2.6$$ 
Note that $B$ is not a Hopf algebra, the coproduct becomes singular at
the points $s=\pm 1.$ 

The total space $P$ is the algebra with generators $(a_{ij},
b_{ij},t)$ where the $a_{ij}$'s commute with the $b_{kl}$'s, the 
$a_{ij}$'s and $b_{ij}$'s satisfy the  $SU_q(N)$ quadratic relations 
(2.1) but the unitarity constraints are 
$$a_{ij}a_{ik}^*= t \delta_{jk}, \,\, b_{ij} b_{ik}^*=(1-t)
\delta_{jk}.\tag2.7$$

$P$ is a $A=SU_q(N)$ comodule algebra. The coaction is 
$$\Delta_R(a_{ij})= a_{ik}\otimes g_{kj}, \,\, \Delta_R(b_{ij}) = b_{ik} 
\otimes g_{kj}, \,\, \Delta_R(t) = t\otimes 1.\tag2.8$$ 

We have an injective algebra homomorphism  $\phi:B\to P$ given by 
$$\phi(g_{ij})= a_{ik}b_{jk}^*, \,\,\, \phi(s)=t.$$ 
This maps to the $SU_q(N)$ coinvariants in $P,$ 
$$\align \Delta(\phi(g_{ij}))&= \Delta_R(a_{ik})\Delta_R(b_{jk}^*)= 
(a_{im}\otimes g_{mk})\cdot (b_{jl}\otimes g_{lk})^* \\
&= a_{im}b_{jl}^*\otimes g_{mk} g_{lk}^* 
=a_{im} b_{jm}^* \otimes 1 .\tag2.9\}\endalign$$ 

It is easy to see that the map $x\otimes y\mapsto (x\otimes
1)\Delta_R(y)$ from $P\otimes P\to P\otimes A$ is surjective 
(freeness condition for $A$ coaction, Axiom 4 in [Br-Ma]) and that 
$B=P^A = \{x\in P| \Delta_R(x) = x\otimes 1\}$ (Axiom 3). 

In addition, we have \it local trivializations \rm $\Phi_{pm}: A\to P_{\pm}$ 
which satisfy  
$$\Delta_R\circ \Phi_{\pm}=(\Phi_{\pm}\otimes id)\circ \Delta$$
and $\Phi_{\pm} (1_A)=1_P.$ In addition, there is a local convolution inverse
$\Phi^{-1}: A\to P,$ that is, 
$$\Phi(a')\Phi^{-1}(a'') = 1_P \text{ for } \Delta(a)=\sum_{(a)}
a'\otimes a''.$$ 
Local means here that we have to localize the algebra $P$ by adding 
a new generator $r_{\pm}$ with $r_+t=1, r_-(1-t)=1.$ 

The maps which are regular at $s=1$ but singular at $s=-1$ 
are then given as 
$$\Phi_+(g_{ij})= a_{ij} \text{ and } \Phi_+^{-1}(g_{ij}) = r_+ a_{ij}.$$ 

There is also a pair of maps $\Psi_-,\Psi_-^{-1}:A\to P_-$ which are regular
at $s=-1$ but singular at $s=1.$ They are given as
$$\Psi_+(g_{ij})= a_{ij}, \,\, \Psi^{-1}(g_{ij}) = r_- b_{ij}.$$

\enddocument